\documentclass[11pt,a4paper]{article}
\usepackage{amsthm,amsmath,amssymb,amsfonts,bbm}
\usepackage{aliascnt}
\usepackage[hidelinks]{hyperref}
\usepackage{geometry}
\usepackage[english]{babel}
\usepackage{microtype}
\usepackage[capitalise]{cleveref}

\newtheorem{theorem}{Theorem}[section]
\newaliascnt{prop}{theorem}
\newtheorem{prop}[prop]{Proposition}
\aliascntresetthe{prop}
\newaliascnt{lemma}{theorem}
\newtheorem{lemma}[lemma]{Lemma}
\aliascntresetthe{lemma}
\newaliascnt{cor}{theorem}
\newtheorem{cor}[cor]{Corollary}
\aliascntresetthe{cor}
\newaliascnt{conj}{theorem}

\aliascntresetthe{conj}
\newaliascnt{prob}{theorem}

\aliascntresetthe{prob}
\newaliascnt{claim}{theorem}

\aliascntresetthe{claim}

\newaliascnt{fact}{theorem}
\newtheorem{fact}[fact]{Fact}
\aliascntresetthe{fact}

\crefname{theorem}{Theorem}{Theorems}
\Crefname{theorem}{Theorem}{Theorems}
\crefname{prop}{Proposition}{Propositions}
\Crefname{prop}{Proposition}{Propositions}
\crefname{lemma}{Lemma}{Lemmas}
\Crefname{lemma}{Lemma}{Lemmas}
\crefname{cor}{Corollary}{Corollaries}
\Crefname{cor}{Corollary}{Corollaries}
\crefname{conj}{Conjecture}{Conjectures}
\Crefname{conj}{Conjecture}{Conjectures}
\crefname{prob}{Problem}{Problems}
\Crefname{prob}{Problem}{Problems}
\crefname{claim}{Claim}{Claims}
\Crefname{claim}{Claim}{Claims}
\crefname{fact}{Fact}{Facts}
\Crefname{fact}{Fact}{Facts}

\theoremstyle{definition}
\newaliascnt{defn}{theorem}
\newtheorem{defn}[defn]{Definition}
\aliascntresetthe{defn}
\crefname{defn}{Definition}{Definitions}
\Crefname{defn}{Definition}{Definitions}
\newaliascnt{exam}{theorem}

\aliascntresetthe{exam}
\newaliascnt{setup}{theorem}

\aliascntresetthe{setup}
\newaliascnt{exer}{theorem}

\aliascntresetthe{exer}
\newaliascnt{ques}{theorem}
\newtheorem{ques}[ques]{Question}
\aliascntresetthe{ques}

\theoremstyle{remark}
\newtheorem*{remark}{Remark}

\crefname{ques}{Question}{Questions}
\Crefname{ques}{Question}{Questions}

\newcommand{\E}{\mathbb{E}}
\newcommand{\one}{\mathbf{1}}
\newcommand{\Cov}{\mathrm{Cov}}
\newcommand{\Var}{\mathrm{Var}}
\newcommand{\Inf}{\mathrm{Inf}}
\newcommand{\I}{\mathcal{I}}
\newcommand{\N}{\mathcal{L}}
\newcommand\tup[1]{\left\langle #1 \right\rangle}

\title{Fourier level inequalities on the discrete cube via complex hypercontractivity and Riesz transforms}
\author{%
Fan Chang\thanks{%
School of Statistics and Data Science, Nankai University, Tianjin, China.
Email: \texttt{1120230060@mail.nankai.edu.cn}.
Supported by the National Natural Science Foundation of China (NSFC)
under grant 124B2019.},\quad
Zhenyu Li\thanks{%
School of Mathematics, Shandong University, Jinan, China.
Email: \texttt{zy.li@mail.sdu.edu.cn}.
},\quad
Guanghui Wang\thanks{%
School of Mathematics, Shandong University, Jinan, China. Supported by National Key R\&D Program of China (2020YFA0712400) and National Natural Science Foundation of China (12231018) Email: \texttt{ghwang@sdu.edu.cn}.
},\quad
Mingrui Wang\thanks{%
School of Mathematics, Shandong University, Jinan, China.
Email: \texttt{202300091103@mail.sdu.edu.cn}.
}.\quad
}
\date{}

\hypersetup{
  pdftitle={A Direct Fourier-Analytic Boolean Level-1:3 Inequality and Antipodal Correlation Bounds},
  pdfauthor={Fan Chang},
  pdfsubject={Discrete Fourier analysis and correlation inequalities on the Boolean hypercube},
  pdfkeywords={Boolean functions, correlation inequalities, complex hypercontractivity, random restrictions}
}

\begin{document}
\maketitle

\begin{abstract}
Eldan's Level-$1{:}3$ inequality is a central tool in his second-order correlation inequalities for increasing Boolean functions. His proof is going to Gaussian space and uses an entropy estimate for second Hermite tensors. We give a direct proof on the discrete hypercube. For increasing real-valued functions $f,g:\{0,1\}^n\to\mathbb R$, we prove
$$
W_3[f,g]\ge -2W_1[f,g]\log\!\left(
\frac{\sqrt{\Var(f)\Var(g)}}{W_1[f,g]}
\right).
$$
The proof combines complex hypercontractivity with the discrete Riesz transforms. We also give a short proof of a level-$1:d$ inequality obtained by Keller and Kindler.

As an application, we recover Eldan's correlation inequality with a slight improvement. If $f,g:\{0,1\}^n\to\{0,1\}$ are increasing and $g$ is antipodal, then
$$
\Cov(f,g)\ge \frac1{24}\,
\frac{\sum_i\Inf_i[f]\Inf_i[g]}{\sqrt{\log\!\left(
4e\sqrt{\Var(f)\Var(g)}/(\sum_i\Inf_i[f]\Inf_i[g])\right)}},
$$
where the square-root logarithm cannot be improved as shown in the OR--majority example.

\vspace{10pt}
\noindent\textbf{Keywords:} Correlation inequalities, Influences, Complex hypercontractivity, Discrete Fourier analysis, Random restrictions

\vspace{10pt}
\noindent\textbf{Mathematics Subject Classification:} 06E30; 42C10; 60C05
\end{abstract}

\section{Introduction}

Level inequalities control the Fourier mass of a Boolean function at a given degree in terms of its expectation, its influences, or lower Fourier levels. A basic example is the Level-$1$ inequality, often called Chang's lemma on the discrete cube~\cite{Chang2002,IMR2014}: for a Boolean function with expectation $\alpha\in(0,1)$ under the uniform measure, the sum of squares of its first-level Fourier coefficients is at most $2\alpha^2\log(1/\alpha)$. More generally, the Level-$d$ inequality bounds the sum of squares of its level-$d$ Fourier coefficients by $\alpha^2(2e\log(1/\alpha)/d)^d$ for integers $1\le d\le2\log(1/\alpha)$. This estimate follows from Bonami--Beckner real hypercontractivity~\cite{Beckner1975,Bonami1970}; see O'Donnell~\cite[Chapter~9]{ODonnell2014} for a systematic survey. Recent work has extended these level inequalities beyond the uniform cube under globalness/pseudorandomness conditions~\cite{KLLM2024,KLM2026}. Green and Sawhney~\cite{GS2025} prove an arithmetic Level-$d$ inequality and use it in additive number theory. 

In probabilistic combinatorics, Talagrand's celebrated Level-$1{:}2$ inequality~\cite{Talagrand1996} and its higher-level extensions~\cite{KK2013} have important applications to correlation inequalities~\cite{KKM2016,Talagrand1996} and noise sensitivity~\cite{BKS1999}. We work throughout on the discrete hypercube $\{0,1\}^n$ equipped with the uniform product measure $\mu=(\frac12\delta_0+\frac12\delta_1)^{\otimes n}$. Unless stated otherwise, all expectations are taken with respect to $\mu$. For functions $f,g:\{0,1\}^n\to\mathbb R$, write
$$
\Cov(f,g):=\mathbb E[fg]-\mathbb E[f]\mathbb E[g],
\qquad
\Var(f):=\Cov(f,f).
$$
We order vectors coordinatewise. A function $f$ is \emph{increasing} if $x\le y$ implies $f(x)\le f(y)$. For $S\subseteq[n]$, let $\chi_S(x):=(-1)^{\sum_{i\in S}x_i}$ and $\hat f(S):=\E[f\chi_S]$. We use the notation
\begin{equation}\label{eq:cross-fourier-level-intro}
W_k[f,g]:=\sum_{|S|=k}\hat f(S)\hat g(S),
\qquad W_k[f]:=W_k[f,f].
\end{equation}

\begin{theorem}[{Talagrand's Level-$1{:}2$ inequality~\cite[Theorem~2.4]{Talagrand1996}}]
There is a universal constant $C>0$ such that for increasing Boolean functions $f,g:\{0,1\}^n\to\{0,1\}$,
\begin{equation}\label{eq:talagrand-level12-intro}
W_2[f,g]\le C W_1[f,g]\log\left(\frac{e}{W_1[f,g]}\right).
\end{equation}
\end{theorem}

One application of \eqref{eq:talagrand-level12-intro} is Talagrand's quantitative form of the Harris--Kleitman inequality~\cite{Harris1960,Kleitman1966}. For $i\in[n]$, the influence of the $i$-th coordinate on a Boolean function $f$ is $\Inf_i[f]:=\mathbb P\bigl[f(x)\ne f(x\oplus e_i)\bigr]$, and we use the cross-total-influence $\I[f,g]:=\sum_{i=1}^n\Inf_i[f]\Inf_i[g]$. Talagrand~\cite{Talagrand1996} proved that
\begin{equation}\label{eq:talagrand-correlation-intro}
\Cov(f,g)\ge C\,\frac{\I[f,g]}{\log(e/\I[f,g])}.
\end{equation}
His proof combines the Level-$1{:}2$ inequality with an induction on the dimension. Kalai, Keller and Mossel~\cite[Section~2.3]{KKM2016} replace this induction by a semigroup argument. Subsequent work developed symmetry, regularity, Gaussian, and second-order refinements; see, for example,~\cite{Eldan2022,FKKK2018,KKM2016,KMS2014}.

Eldan~\cite{Eldan2022} showed that a lower bound between the first and third levels. His result is formulated on Gaussian space. Let $\gamma_n$ be standard Gaussian measure on $\mathbb R^n$, and define the first and third Hermite tensors by
$$
H^{(1)}(x)=x,\qquad
H^{(3)}_{ijk}(x)=x_ix_jx_k-\delta_{ij}x_k-\delta_{ik}x_j-\delta_{jk}x_i,
$$
where $\delta_{ij}$ is the Kronecker delta. For $F\in L^2(\gamma_n)$, let
$$
Q^{(k)}(F):=\int_{\mathbb R^n}F(x)H^{(k)}(x)\,d\gamma_n(x).
$$
\begin{theorem}[Eldan's Gaussian Level-$1{:}3$ inequality {\cite[Proposition~17]{Eldan2022}}]
Let $F,G\in L^2(\gamma_n)$ be coordinatewise increasing. Then
\begin{equation}\label{eq:eldan-gaussian-level13}
\langle Q^{(3)}(F),Q^{(3)}(G)\rangle
\ge -e^8 \langle Q^{(1)}(F),Q^{(1)}(G)\rangle\log\left(\frac{e\sqrt{\Var_{\gamma_n}(F)\Var_{\gamma_n}(G)}}{\langle Q^{(1)}(F),Q^{(1)}(G)\rangle}\right).
\end{equation}
\end{theorem}
Eldan proves this inequality from an entropy bound for second Hermite tensors using Ornstein--Uhlenbeck smoothing and Gaussian integration by parts~\cite[Section~3]{Eldan2022}. To obtain correlation inequalities, Eldan applies this inequality to conditional Gaussian functions along a Brownian martingale and controls the resulting differential inequality. An orthant lift transfers the correlation bounds to the discrete hypercube. In discrete isoperimetry, Eldan and Gross~\cite{EG2022} used a related approach based on pathwise stochastic analysis to prove a Talagrand's isoperimetric conjecture~\cite{Talagrand1997}. Subsequent work has provided more direct proofs that avoid this machinery, using semigroup gradient estimates in van Handel's argument~\cite{Rosenthal2020}, random restrictions and Fourier tail
bounds~\cite{EKLM2025}, or the Bobkov--G\"otze isoperimetric inequality~\cite{IvanisviliZhang2026}. 

Motivated by these works, we prove the following Boolean Level-$1{:}3$ inequality. 
\begin{theorem}[Boolean Level-$1{:}3$ inequality]\label{thm:level13-main}
Let $f,g:\{0,1\}^n\to\mathbb R$ be increasing. If $W_1[f,g]>0$, then
\begin{equation}\label{eq:level13-main}
W_3[f,g]
\ge -2W_1[f,g]\log\left(
\frac{\sqrt{\Var(f)\Var(g)}}{W_1[f,g]}
\right).
\end{equation}
If $W_1[f,g]=0$, then $W_3[f,g]=0$.
\end{theorem}
\begin{remark}
By Parseval and Cauchy--Schwarz, $W_1[f,g]\le\sqrt{\Var(f)\Var(g)}$, so the logarithm is non-negative. For increasing Boolean functions, $W_1[f,g]=\I[f,g]/4$; thus \cref{thm:level13-main} gives
$$
W_3[f,g]\ge -\frac{\I[f,g]}2\log\left(\frac{4\sqrt{\Var(f)\Var(g)}}{\I[f,g]}\right).
$$
A Boolean estimate with an unspecified universal constant follows from Eldan's Gaussian Level-1:3 inequality~\cite[Proposition 17]{Eldan2022} through orthant lifting. Here we prove it directly on the cube with constant $2$.
\end{remark}

Our proof uses complex hypercontractivity and discrete Riesz transforms, without passing to Gaussian space. Both tools are well studied on the
cube~\cite{Beckner1975,IvanisviliNazarov2022,LustPiquard1998,LustPiquard2004}, but, to our knowledge, their joint use in quantitative correlation inequalities on the cube is new. This combination may also be useful in Boolean function analysis and extremal combinatorics.

As a further application of the Riesz-transform method, we give a short proof of a two-function higher-level inequality. Such estimates were proved by Keller and Kindler~\cite[Lemma~5]{KK2013} and are used in~\cite[Lemma~2.8]{KKM2016}.

\begin{theorem}[A Level-$1{:}d$ inequality]\label{thm:keller-kindler-main}
Let $f,g:\{0,1\}^n\to[-1,1]$ be increasing, and let $d\ge2$ be an integer. If $0<W_1[f,g]\le e^{-2(d-1)}$, then
\begin{equation}\label{eq:keller-kindler-main}
|W_d[f,g]|\le\frac{2d-1}{d^2}\left(\frac{2e}{d-1}\right)^{d-1}\,W_1[f,g]\left(\log\frac{d}{W_1[f,g]}\right)^{d-1}.
\end{equation}
If $W_1[f,g]=0$, then $W_d[f,g]=0$.
\end{theorem}

The proof in~\cref{subsec:keller-kindler} uses `shifted' discrete Riesz transforms and real hypercontractivity.

\subsection{Correlation inequalities}

Our main application is a direct proof of Eldan's correlation inequality, with an explicit constant. Recall that a Boolean function $g$ is \emph{antipodal} if $g(x)+g(\one-x)=1$ for every $x\in\{0,1\}^n$. 

\begin{theorem}\label{thm:antipodal-main}
Let $f,g:\{0,1\}^n\to\{0,1\}$ be increasing, and suppose that $g$ is antipodal. Then
\begin{equation}\label{eq:antipodal-main}
\Cov(f,g)\ge \frac1{24}\cdot
\frac{\I[f,g]}{\sqrt{\log\left(
4e\sqrt{\Var(f)\Var(g)}/\I[f,g]
\right)}}.
\end{equation}
The right-hand side is interpreted as $0$ when $\I[f,g]=0$.
\end{theorem}

Since an antipodal Boolean function has expectation $\frac{1}{2}$ and variance $\frac{1}{4}$, \eqref{eq:antipodal-main} may also be written as
$$
\Cov(f,g)\ge \frac1{24}\cdot
\frac{\I[f,g]}{\sqrt{\log\left(2e\sqrt{\Var(f)}/\I[f,g]\right)}}.
$$
Moreover, $\sqrt{\Var(f)\Var(g)}\le\frac{1}{4}$, and hence \cref{thm:antipodal-main} implies Eldan's result~\cite[Theorem 2]{Eldan2022} with an explicit constant:
\begin{equation}\label{eq:eldan-explicit}
\Cov(f,g)
\ge \frac1{24}\cdot
\frac{\I[f,g]}{\sqrt{\log\left(2e/\I[f,g]\right)}}.
\end{equation}

The same method retains the sign of the second Fourier level. The second Fourier level contributes to the lower bound with its original sign.

\begin{theorem}[One-sided second-order refinement]\label{thm:second-order-main}
Let $f,g:\{0,1\}^n\to\{0,1\}$ be increasing. If $W_1[f,g]>0$, then
\begin{equation}\label{eq:second-order-main}
\Cov(f,g)
\ge \frac{W_1[f,g]}{16}
\min\left\{
\frac1{\sqrt{\log\left(
 e\sqrt{\Var(f)\Var(g)}/W_1[f,g]
\right)}},
\frac{W_1[f,g]}{\bigl(-W_2[f,g]\bigr)_+}
\right\},
\end{equation}
where the second term is interpreted as $+\infty$ when $W_2[f,g]\ge0$. If $W_1[f,g]=0$, then $\Cov(f,g)\ge0$.
\end{theorem}

For comparison, after changing Fourier normalization, Eldan's stated second-order theorem~\cite[Theorem~3]{Eldan2022} gives
$$
\Cov(f,g)\ge CW_1[f,g]\min\left\{
\frac1{\sqrt{\log(e/W_1[f,g])}},\frac{W_1[f,g]}{|W_2[f,g]|}\right\}.
$$

\begin{remark}[Chv\'atal's conjecture and the dream inequality]
The motivation for Eldan's results comes from Chv\'atal's conjecture: every hereditary family has a largest intersecting subfamily consisting of all its members containing a fixed element. Friedgut, Kahn, Kalai and Keller~\cite[Proposition~2.1]{FKKK2018} showed that the conjecture is equivalent to
\begin{equation}\label{eq:chvatal-correlation-intro}
\Cov(f,g)\ge\frac14\min_i\Inf_i[f]
\end{equation}
for increasing Boolean $f,g$ with $g$ antipodal. The conjecture has recently been proved by Chang, Liu and Liu~\cite{CLL2026}, through a stronger spectral correlation inequality.

A natural, stronger candidate was the ``dream inequality''
\begin{equation}\label{eq:dream-intro}
\Cov(f,g)\ge\frac14\I[f,g].
\end{equation}
Indeed, a balanced Boolean function satisfies $\sum_i\Inf_i[g]\ge1$ by Harper's edge-isoperimetric inequality, so \eqref{eq:dream-intro} would imply \eqref{eq:chvatal-correlation-intro}; this implication is recorded in~\cite[Proposition~2.2]{FKKK2018}. The dream inequality is false, even when $g$ is antipodal~\cite{Chang2026}. Indeed, the square-root logarithm has the correct order even in the antipodal setting.
\end{remark}

\begin{prop}[Sharpness of the logarithmic exponent]\label{prop:sharpness}
For every $\alpha<1/2$, there is no universal constant $c>0$ such that
$$
\Cov(f,g)\ge c\cdot\frac{\I[f,g]}{\log^{\alpha}(2e/\I[f,g])}
$$
for all increasing Boolean $f,g$ with $g$ antipodal.
\end{prop}

\begin{remark}
De, Nadimpalli and Servedio~\cite{DNS2022} also use complex analysis to prove quantitative correlation inequalities. They apply Hadamard's three-circles theorem to obtain an extremal bound for power series
with bounded sums of absolute coefficients, and then use monotonicity
along the noise semigroup to deduce a correlation bound. On the uniform cube their argument gives a denominator of order $\log^2(e/W_1[f,g])$. Our use of
complex noise is different. We use Beckner’s complex hypercontractive inequality only in the proof of~\cref{lem:imaginary-level-two}.  
\end{remark}
\subsection{Proof overview}
The proof of~\cref{thm:level13-main} has two ingredients. 

\medskip\noindent\emph{Imaginary noise.}
Let $f,g\ge0$ satisfy $\E f=\E g=1$, and extend $T_z$ to
$z\in\mathbb C$ by its Fourier formula. As $r\to0$,
$$
\bigl|\E[fT_{ir}g]\bigr|^2=1+\bigl(W_1[f,g]^2-2W_2[f,g]\bigr)r^2+O(r^4).
$$
For $0<r<1$, Beckner's complex hypercontractive~\cite{Beckner1975} gives $\|T_{ir}g\|_{p'}\le\|g\|_p$, where $p=1+r^2$ and
$p'=\frac{p}{p-1}$. H\"older's inequality and $L^1$--$L^2$
interpolation therefore give
$$
\bigl|\E[fT_{ir}g]\bigr|^2\le\bigl(\|f\|_2\|g\|_2\bigr)^{\frac{4r^2}{1+r^2}}.
$$
Comparing quadratic terms yields
$$
W_2[f,g]\ge-2\log(\|f\|_2\|g\|_2)+\frac12W_1[f,g]^2.
$$

\medskip\noindent\emph{Discrete Riesz transforms.}
For increasing $f,g$, we apply the homogeneous form of the
preceding estimate to
$$
\mathfrak R_i f:=\frac12\partial_i\N^{-1/2}(f-\E f),
\qquad \N\chi_S=|S|\chi_S,
$$
where
$\partial_i f=f(x^{(i\to1)})-f(x^{(i\to0)})$.
These oriented discrete Riesz transforms are non-negative
and satisfy
$$
\begin{aligned}
\sum_i\E[\mathfrak R_i f]\E[\mathfrak R_i g]&=W_1[f,g],\\
\sum_i\|\mathfrak R_i f\|_2^2&=\Var(f),\\
\sum_iW_2[\mathfrak R_i f,\mathfrak R_i g]&=W_3[f,g].
\end{aligned}
$$
Note that the usual half-derivatives give
$\sum_{S\ne\emptyset}|S|\hat f(S)^2$ instead of $\Var(f)$ in the second identity. The multiplier $|S|^{-1/2}$ cancels this multiplicity exactly. Jensen's inequality and Cauchy--Schwarz then yield~\cref{thm:level13-main}.

For the correlation inequalities, let $f,g$ be increasing Boolean
functions. Motivated by the proof of Theorem~1.3 in~\cite{KKM2016},
we estimate $\Cov(f,T_\rho g)$ at a suitable small parameter and then use
$$
\Cov(f,g)\ge\Cov(f,T_\rho g)=\sum_{k\ge1}\rho^kW_k[f,g],
$$
where the inequality follows from~\cite[Lemma~1.9]{KKM2016}. Using Taylor's formula, \cref{lem:beta-extraction} gives
\begin{equation}\label{eq:decomp-identity}
\Cov(f,T_\rho g)=\rho W_1[f,g]+\rho^2W_2[f,g]+3\int_0^\rho\frac{(\rho-p)^2}{(1-p)^3}
\E_{(J,z)}W_3[f_{J\to z},g_{J\to z}]\,dp,
\end{equation}
where $f_{J\to z}$ and $g_{J\to z}$ are obtained by applying the same random restriction: each coordinate belongs to $J$ independently with probability $p$, and the coordinates in $J$ are fixed to a uniform assignment $z\in\{0,1\}^J$. Thus $f_{J\to z}(y)=f(x_J=z,x_{J^c}=y)$ for
$y\in\{0,1\}^{J^c}$, and similarly for $g$. A Fourier calculation analogous to that in the proof of~\cite[Fact~2.4]{EKLM2025} relates the expected cross levels of the restrictions to the original functions: for every integer $r\ge1$,
\begin{equation}\label{eq:restriction-overview}
\E_{(J,z)}W_r[f_{J\to z},g_{J\to z}]=\frac{(1-p)^r}{r!}\left.\frac{d^r}{d\rho^r}\Cov(f,T_\rho g)\right|_{\rho=p}.
\end{equation}
We then apply \cref{thm:level13-main} inside the integral of~\eqref{eq:decomp-identity} and use the fact that the expected variance of a restriction is at
most the original variance. It gives
$$
\begin{aligned}
\Cov(f,T_\rho g)&\ge \rho W_1[f,g]+\rho^2W_2[f,g]\\
&\quad-6\rho^2\Cov(f,T_\rho g)
\log\left(\frac{\sqrt{\Var(f)\Var(g)}\log(1/(1-\rho))}{\Cov(f,T_\rho g)}\right).
\end{aligned}
$$
When $g$ is antipodal, $W_2[f,g]=0$. Taking a suitable $\rho$ bounds $\Cov(f,T_\rho g)$ below by a fixed positive multiple of $\rho W_1[f,g]$, proving~\cref{thm:antipodal-main}.

\medskip\noindent\emph{Organization.}
\Cref{sec:prelim} reviews Fourier analysis on the discrete cube, complex hypercontractive and random restrictions. In \cref{sec:imaginary} we prove the imaginary-noise estimate. In \cref{sec:fractional-gradient} we apply it to discrete Riesz transforms to prove~\cref{thm:level13-main}, and use shifted Riesz transforms to prove~\cref{thm:keller-kindler-main}. \Cref{sec:beta} gives the proofs of correlation bounds. We prove \cref{prop:sharpness} and discuss the optimal constant in the level inequality in \cref{sec:sharpness}.

\section{Preliminaries}\label{sec:prelim}

\subsection{Fourier analysis on the discrete cube}

We equip $L^2(\{0,1\}^n)$ with the inner product $\tup{f,g}:=\mathbb E[fg]$ and let $\Vert f\Vert_q:= \left ( \mathbb E |f|^q \right )^{1/q}$ for $1\le q<\infty$. The characters $\{\chi_S\}_{S\subseteq[n]}$ form an orthonormal basis, and
$$
f=\sum_{S\subseteq[n]}\hat f(S)\chi_S,
\qquad
\hat f(S)=\tup{f,\chi_S}.
$$
For $k\ge0$, the $k$-th cross Fourier level is
$$
W_k[f,g]:=\sum_{|S|=k}\hat f(S)\hat g(S).
$$
Parseval gives $\Var(f)=\sum_{S\ne\emptyset}\hat f(S)^2$. Consequently, by Cauchy--Schwarz,
$$
\begin{aligned}
|W_1[f,g]|
&\le\left(\sum_{i=1}^n\hat f(\{i\})^2\right)^{1/2}\left(\sum_{i=1}^n\hat g(\{i\})^2\right)^{1/2}\\
&\le\left(\sum_{S\ne\emptyset}\hat f(S)^2\right)^{1/2}\left(\sum_{S\ne\emptyset}\hat g(S)^2\right)^{1/2}=\sqrt{\Var(f)\Var(g)}.
\end{aligned}
$$

For $i\in[n]$, define the discrete derivative
$$
\partial_i f(x):=f(x^{(i\to1)})-f(x^{(i\to0)}),
\qquad
x\in\{0,1\}^{[n]\setminus\{i\}}.
$$
Then $f$ is increasing if and only if $\partial_i f\ge0$ pointwise for every $i$. A direct Fourier calculation gives
\begin{equation}\label{eq:derivative-fourier}
\partial_i f
=-2\sum_{S\ni i}\hat f(S)\chi_{S\setminus\{i\}}.
\end{equation}
In particular,
\begin{equation}\label{eq:mean-derivative}
\mathbb E[\partial_i f]=-2\hat f(\{i\}).
\end{equation}
For a Boolean function $f$, the influence of coordinate $i$ is
$$
\Inf_i[f]
:=\mathbb P\bigl[f(x)\ne f(x\oplus e_i)\bigr]
=\mathbb E[|\partial_i f|^2].
$$
If $f$ is also increasing, then $\partial_i f\in\{0,1\}$, so
\begin{equation}\label{eq:influence-fourier-one}
\Inf_i[f]=\mathbb E[\partial_i f]=-2\hat f(\{i\}).
\end{equation}
Consequently, for increasing Boolean $f,g$,
\begin{equation}\label{eq:W1-influence}
W_1[f,g]=\frac14\I[f,g].
\end{equation}

\subsection{Noise operators and the number operator}

For $z\in\mathbb C$, define the complex noise operator by
$$
T_z f:=\sum_{S\subseteq[n]}z^{|S|}\hat f(S)\chi_S.
$$
For $0\le\rho\le1$, $T_\rho$ is a positive Markov operator. We also write $P_t:=T_{e^{-t}}$ for $t\ge0$. The number operator $\N$ is defined by $\N\chi_S=|S|\chi_S$, so $P_t=e^{-t\N}$. On mean-zero functions, define
$$
\N^{-1/2}f:=\sum_{S\ne\emptyset}\frac{\hat f(S)}{\sqrt{|S|}}\chi_S.
$$
\begin{fact}
For a real-valued function $f:\{0,1\}^n\to\mathbb{R}$ with $\E[f]=0$,
\begin{equation}\label{eq:fractional-integral}
\N^{-1/2}f
=\frac1{\sqrt\pi}\int_0^\infty t^{-1/2}P_t f\,dt.
\end{equation}
\end{fact}
\begin{proof}
Since $\Gamma(\frac{1}{2})=\int_0^\infty u^{-\frac{1}{2}}e^{-u}\,du=\sqrt\pi$, the substitution $u=\lambda t$ gives, for $\lambda>0$,
$$
\frac1{\sqrt\pi}\int_0^\infty t^{-\frac{1}{2}}e^{-\lambda t}\,dt
=\frac1{\sqrt{\pi\lambda}}\int_0^\infty u^{-\frac{1}{2}}e^{-u}\,du
=\lambda^{-\frac{1}{2}}.
$$
Applying this identity to each nonconstant Fourier character gives~\eqref{eq:fractional-integral}.
\end{proof}
\begin{remark}
Indeed, we demand the integral $\int_0^\infty t^{-1/2}\|P_tf\|_2\,dt$ converges in our proof. By definition, we obtain that 
$$\|P_tf\|_2^2=\sum_{S\subseteq [n]}e^{-2t|S|}\hat f(S)^2 \le  e^{-2t}\|f\|_2^2 + (1-e^{-2t}) \hat f(\emptyset)^2 .$$
Therefore,
$$ \int_0^\infty t^{-1/2}\|P_tf\|_2\,dt\le\sqrt\pi\,\|f\|_2 + \hat f(\emptyset)^2 \int_0^\infty (1-e^{-2t}) t^{-1/2} dt.
$$
Observe that $\int_0^\infty (1-e^{-2t}) t^{-1/2} dt$ diverges. However, when the function $f$ is mean-zero, that is, $\mathbb E[f] = \hat f(\emptyset) =0$, the second term cancels out, and the desired integral converges.

\end{remark}
The derivative and heat semigroup satisfy
\begin{equation}\label{eq:derivative-commutation}
\partial_iP_t f=e^{-t}P_t^{(-i)}(\partial_i f),
\end{equation}
where $P_t^{(-i)}$ is the heat semigroup on $\{0,1\}^{[n]\setminus\{i\}}$. To verify the factor $e^{-t}$, use \eqref{eq:derivative-fourier}:
$$
\partial_iP_tf=-2\sum_{S\ni i}e^{-t|S|}\hat f(S)\chi_{S\setminus\{i\}}=e^{-t}\left(-2\sum_{S\ni i}e^{-t(|S|-1)}\hat f(S)\chi_{S\setminus\{i\}}\right)
=e^{-t}P_t^{(-i)}(\partial_i f).
$$

We shall use the following endpoint of complex hypercontractivity.

\begin{theorem}[Beckner's complex hypercontractive endpoint {\cite[Lemma~1]{Beckner1975}}]\label{thm:complex-hc}
Let $1<p\le2$ and $p'=p/(p-1)$. For every $f:\{0,1\}^n\to\mathbb C$,
\begin{equation}\label{eq:complex-hc}
\left\|T_{i\sqrt{p-1}}f\right\|_{p'}\le\|f\|_p.
\end{equation}
\end{theorem}

Beckner proves the corresponding one-coordinate two-point inequality; tensorization gives \eqref{eq:complex-hc} in every dimension. For a modern account of this endpoint, its place in the Weissler region, and the tensor-power reduction to the two-point inequality, see \cite[Section~1.2 and Lemma~2.1]{IvanisviliNazarov2022}. We only use the limit $p\downarrow1$.

For real-valued $f,g$ and $0\le\rho\le1$, define
\begin{equation}\label{eq:H-definition}
h(\rho):=\Cov(f,T_\rho g)
=\sum_{k=1}^n\rho^kW_k[f,g].
\end{equation}
For increasing $f,g$, the following monotonicity is due to Kalai, Keller and Mossel~\cite[Lemma~1.9 and Remark~2.5]{KKM2016}. We include the Fourier proof in our normalization. By~\eqref{eq:derivative-fourier},
\begin{equation}\label{eq:H-monotone}
h'(\rho)=\sum_{S\ne\emptyset}|S|\rho^{|S|-1}\hat f(S)\hat g(S)=\frac14\sum_{i=1}^n\tup{\partial_i f,T_\rho\partial_i g}\ge0.
\end{equation}
Here $T_\rho$ in the $i$-th summand acts on the cube without coordinate $i$. The last inequality follows from $\partial_i f,\partial_i g\ge0$ and positivity of $T_\rho$. Hence
\begin{equation}\label{eq:cov-Hdelta}
\Cov(f,g)=h(1)\ge h(\rho).
\end{equation}

\subsection{Antipodality and random restrictions}

If $g$ is antipodal, then $\mathbb Eg=1/2$. Moreover, for every nonempty $S$,
$$
\hat g(S)
=\mathbb E[g(\one-x)\chi_S(\one-x)]
=-(-1)^{|S|}\hat g(S).
$$
Hence
\begin{equation}\label{eq:antipodal-even-vanish}
\hat g(S)=0~\text{whenever }|S|\text{ is positive and even}.
\end{equation}
In particular, $W_2[f,g]=0$ for every $f$.

\medskip
\noindent One of the most useful tools in the analysis of Boolean functions is the notion of restrictions and random restrictions. Let $J\subseteq[n]$ be a set of coordinates thought of as ``restricted'', and let $J^c=[n]\setminus J$ be the set of ``free'' coordinates.

\begin{defn}[Restrictions]\label{def:restriction}
Let $f:\{0,1\}^n\to\mathbb R$, $J\subseteq[n]$, and $z\in\{0,1\}^J$. We denote the restriction by $(J,z)$. The \emph{restricted function} $f_{J\to z}:\{0,1\}^{J^c}\to\mathbb R$ is defined by
$$
f_{J\to z}(y):=f(x_J=z,x_{J^c}=y)=f(y,z),
\qquad y\in\{0,1\}^{J^c}.
$$
In $f(y,z)$, the vector $y$ occupies the coordinates in $J^c$ and $z$ occupies those in $J$. 
\end{defn}

\noindent First, we give a formula for the Fourier coefficients of a restricted function.

\begin{prop}\label{prop:RFour1}
Let $f:\{0,1\}^n\to\mathbb R$, $J\subseteq[n]$, $z\in\{0,1\}^J$, and $T\subseteq J^c$. Then
$$
\widehat{f_{J\to z}}(T)
=\sum_{S\subseteq J}\hat f(S\cup T)\chi_S(z).
$$
\end{prop}

\begin{proof}
Decomposing each Fourier character into its $J$ and $J^c$ parts gives
$$
f(x)
=\sum_{T\subseteq J^c}\sum_{S\subseteq J}
\hat f(S\cup T)\chi_S(x_J)\chi_T(x_{J^c}).
$$
Substitute $x_J=z$ and $x_{J^c}=y$ to obtain
$$
f_{J\to z}(y)
=\sum_{T\subseteq J^c}
\left(\sum_{S\subseteq J}\hat f(S\cup T)\chi_S(z)\right)\chi_T(y).
$$
The coefficient of $\chi_T(y)$ is $\widehat{f_{J\to z}}(T)$ by uniqueness of the Fourier expansion on $\{0,1\}^{J^c}$.
\end{proof}

\begin{defn}[$p$-random restrictions]\label{def:random-restriction}
Let $p\in[0,1]$. A \emph{$p$-random restriction} of $f:\{0,1\}^n\to\mathbb R$ is obtained by independently fixing each coordinate with probability $p$, assigning every fixed coordinate an independent uniform value in $\{0,1\}$, and leaving all remaining coordinates free.
\end{defn}

We write $J\sim\mu_p$ for the random subset of $[n]$ in which each coordinate belongs to $J$ independently with probability $p$; conditional on $J$, we write $z\sim\{0,1\}^J$ for the uniform assignment. Thus a given pair $(J,z)$ has probability $p^{|J|}(1-p)^{n-|J|}2^{-|J|}$.

Restrictions preserve monotonicity and Boolean-valuedness, but need not preserve antipodality. We shall use the following property, which is arbitrary-level two-function version of the calculation in~\cite[Fact~2.4]{EKLM2025}.

\begin{lemma}\label{lem:restriction-level}
Let $(J,z)$ be a $p$-random restriction, with $0\le p<1$. For every integer $r\ge1$,
\begin{equation}\label{eq:restriction-level}
\E_{(J,z)} W_r \left [f_{J\to z},g_{J\to z} \right]
=\frac{(1-p)^r}{r!}h^{(r)}(p).
\end{equation}
\end{lemma}

\begin{proof}
Condition on $J$. By \cref{prop:RFour1}, for $T\subseteq J^c$,
$$
\E_z\left[\widehat{f_{J\to z}}(T)\widehat{g_{J\to z}}(T)\right]=\sum_{S,S'\subseteq J}\hat f(S\cup T)\hat g(S'\cup T)\E_z[\chi_S(z)\chi_{S'}(z)]=\sum_{S\subseteq J}\hat f(S\cup T)\hat g(S\cup T).
$$
The second equality uses orthogonality under the uniform distribution of $z$: the expectation is $1$ for $S=S'$ and $0$ otherwise. Summing over $T\subseteq J^c$ with $|T|=r$ gives
$$
\E_z W_r[f_{J\to z},g_{J\to z}] =\sum_{\substack{T\subseteq J^c\\|T|=r}}\sum_{S\subseteq J}\hat f(S\cup T)\hat g(S\cup T)=\sum_{R\subseteq[n]}\mathbbm{1}_{\{|R\cap J^c|=r\}}\hat f(R)\hat g(R).
$$
Indeed, each $R$ in the last sum is uniquely decomposed as $R=(R\cap J)\cup(R\cap J^c)$.

For a fixed $R$ with $|R|=k\ge r$, the event $|R\cap J^c|=r$ means that exactly $r$ coordinates of $R$ remain free. There are $\binom{k}{r}$ choices for them; each choice has probability $(1-p)^rp^{k-r}$. Coordinates outside $R$ are unrestricted in this event. Therefore
$$
\mathbb P_J(|R\cap J^c|=r)=\binom{k}{r}(1-p)^rp^{k-r},
$$
and averaging over $J$ yields
$$
\E_{(J,z)}W_r[f_{J\to z},g_{J\to z}]=\sum_{R\subseteq[n]}\mathbb P_J(|R\cap J^c|=r)\hat f(R)\hat g(R)=(1-p)^r\sum_{k=r}^n\binom{k}{r}p^{k-r}W_k[f,g].
$$
Since
$$
h^{(r)}(p)
=\sum_{k=r}^n\frac{k!}{(k-r)!}p^{k-r}W_k[f,g]
=r!\sum_{k=r}^n\binom{k}{r}p^{k-r}W_k[f,g],
$$
this proves \eqref{eq:restriction-level}. For $r>n$, both sides vanish.
\end{proof}

The next bound is also an immediate consequence of the law of total variance, first conditioning on the coordinates in a fixed $J$ and then averaging over $J$. We give the Fourier proof as well.

\begin{lemma}\label{lem:restriction-variance}
For a $p$-random restriction,
$$
\E_{(J,z)}\Var(f_{J\to z})\le\Var(f),
\qquad
\E_{(J,z)}\Var(g_{J\to z})\le\Var(g).
$$
\end{lemma}

\begin{proof}
For $p=1$, the restricted functions are constant. For $0\le p<1$, Parseval on the restricted cube and~\cref{lem:restriction-level} give
$$
\begin{aligned}
\E_{(J,z)}\Var(f_{J\to z})
&=\sum_{r=1}^n\E_{(J,z)}W_r[f_{J\to z},f_{J\to z}]=\sum_{r=1}^n\sum_{k=r}^n\binom{k}{r}(1-p)^rp^{k-r}W_k[f,f]\\
&=\sum_{k=1}^n\left(\sum_{r=1}^k\binom{k}{r}(1-p)^rp^{k-r}\right)W_k[f,f]=\sum_{k=1}^n\left(1-p^k\right)W_k[f,f]\\
&=\sum_{S\ne\emptyset}(1-p^{|S|})\hat f(S)^2\le\sum_{S\ne\emptyset}\hat f(S)^2=\Var(f),
\end{aligned}
$$
where we use the binomial theorem in the fourth equality.
\end{proof}

\section{Imaginary noise and a bilinear \texorpdfstring{Level-$2$}{Level-2} inequality}\label{sec:imaginary}

At an imaginary noise parameter, the first Fourier level contributes an imaginary linear term and the second contributes a negative real quadratic term. Complex hypercontractivity bounds the modulus of this expansion. This is the discrete counterpart of the lower bound for second Hermite tensors in~\cite[Lemma~13]{Eldan2022}. 

\begin{lemma}\label{lem:imaginary-level-two}
Let $f,g:\{0,1\}^n\to[0,\infty)$ satisfy $\E f=\E g=1$. Then
\begin{equation}\label{eq:imaginary-level-two}
W_2[f,g]\ge-2\log\bigl(\|f\|_2\|g\|_2\bigr)+\frac12W_1[f,g]^2.
\end{equation}
In particular,
$$
W_2[f,g]\ge-2\log\bigl(\|f\|_2\|g\|_2\bigr).
$$
\end{lemma}

\begin{proof}
For $z\in\mathbb C$, define
$$
F(z):=\langle f,T_zg \rangle =\sum_{d=0}^nW_d[f,g]z^d.
$$
Since $W_0[f,g]=\E f\,\E g=1$, as $r\to0$,
$$
F(ir)=1+iW_1[f,g]r-W_2[f,g]r^2-iW_3[f,g]r^3+O(r^4).
$$
All coefficients $W_d[f,g]$ are real. Squaring the real and imaginary parts separately, we obtain
\begin{equation}\label{eq:F-expansion}
\begin{aligned}
|F(ir)|^2
&=\bigl(1-W_2[f,g]r^2+O(r^4)\bigr)^2+\bigl(W_1[f,g]r-W_3[f,g]r^3+O(r^4)\bigr)^2\\
&=1+\bigl(W_1[f,g]^2-2W_2[f,g]\bigr)r^2+O(r^4).
\end{aligned}
\end{equation}

Next, we give an upper bound for $|F(ir)|$. For $0<r<1$, take $p=1+r^2$ and $p'=\frac{p}{p-1}$. H\"older's inequality and \cref{thm:complex-hc} give
$$
|F(ir)|\le\|f\|_p\|T_{ir}g\|_{p'}\le\|f\|_p\|g\|_p.
$$
For $1<p<2$, another application of H\"older, with exponents $\frac{1}{2-p}$ and $\frac{1}{p-1}$, gives
$$
\E|f|^p=\E\bigl[|f|^{2-p}(|f|^2)^{p-1}\bigr]\le(\E|f|)^{2-p}(\E|f|^2)^{p-1}.
$$
Since $f\ge0$ and $\E f=1$, it follows that
$$
\|f\|_p\le\|f\|_1^{\frac{2-p}{p}}\|f\|_2^{\frac{2(p-1)}{p}}=\|f\|_2^{\frac{2r^2}{1+r^2}}.
$$
Applying the same estimate to $g$, we obtain
\begin{equation}\label{eq:F-upper}
    \begin{split}
|F(ir)|^2&\le\bigl(\|f\|_2\|g\|_2\bigr)^{\frac{4r^2}{1+r^2}}=\exp\left(\frac{4r^2}{1+r^2}\log\bigl(\|f\|_2\|g\|_2\bigr)\right)\\
&=\exp\left(4r^2\log\bigl(\|f\|_2\|g\|_2\bigr)+O(r^4)\right)=1+4r^2\log\bigl(\|f\|_2\|g\|_2\bigr)+O(r^4),
\end{split}
\end{equation}
where we use the fact that $\frac{4r^2}{1+r^2}=4r^2+O(r^4)$ and the expansion $e^v=1+v+O(v^2)$. Comparing~\eqref{eq:F-expansion} with~\eqref{eq:F-upper}, and letting $r\downarrow0$ give
$$
W_1[f,g]^2-2W_2[f,g]\le4\log\bigl(\|f\|_2\|g\|_2\bigr),
$$
which is \eqref{eq:imaginary-level-two}.
\end{proof}

We shall use the following homogeneous form.

\begin{cor}\label{cor:homogeneous-level-two}
Let $f,g:\{0,1\}^n\to[0,\infty)$. If $\E f\,\E g>0$, then
\begin{equation}\label{eq:homogeneous-level-two}
W_2[f,g]\ge -2\E f\,\E g\log\left(\frac{\|f\|_2\|g\|_2}{\E f\,\E g}\right)+\frac{W_1[f,g]^2}{2\E f\,\E g}.
\end{equation}
In particular,
\begin{equation}\label{eq:homogeneous-level-two-simple}
W_2[f,g]\ge -2\E f\,\E g\log\left(
\frac{\|f\|_2\|g\|_2}{\E f\,\E g}
\right).
\end{equation}
If $\E f=0$ or $\E g=0$, then one of $f,g$ vanishes identically.
\end{cor}

\begin{proof}
If $\E f\,\E g>0$, apply \cref{lem:imaginary-level-two} to $f/\E f$ and $g/\E g$, and multiply the resulting inequality by $\E f\,\E g$. If $\E f=0$, then the non-negative function $f$ is identically zero; the case $\E g=0$ is identical.
\end{proof}

\section{Discrete Riesz transforms and the \texorpdfstring{Level-$1{:}3$}{Level-1:3} inequality}\label{sec:fractional-gradient}

To obtain a Level-1:3 inequality from \eqref{eq:homogeneous-level-two-simple}, a direct thought is to use derivatives. For increasing $f$, the ordinary half-derivatives $\frac12\partial_i f$ are non-negative. By \eqref{eq:derivative-fourier} and Parseval,
$$
\begin{aligned}
\sum_{i=1}^n\left\|\frac12\partial_i f\right\|_2^2
&=\sum_{i=1}^n\sum_{S\ni i}\hat f(S)^2
=\sum_{S\ne\emptyset}|S|\hat f(S)^2,\\
\sum_{i=1}^nW_2\left[\frac12\partial_i f,\frac12\partial_i g\right]
&=\sum_{i=1}^n\sum_{\substack{S\ni i\\|S|=3}}\hat f(S)\hat g(S)
=3W_3[f,g].
\end{aligned}
$$
The first identity counts each $S$ once for every $i\in S$, and the second counts each three-element set three times. Also,
$\E[\frac12\partial_i f]\E[\frac12\partial_i g]=\hat f(\{i\})\hat g(\{i\})$.

Suppose $W_1[f,g]>0$. Applying \eqref{eq:homogeneous-level-two-simple} to the latter equation, then using Jensen with weights $\frac{\hat f(\{i\})\hat g(\{i\})}{W_1[f,g]}$, gives
$$
3W_3[f,g]\ge-2\sum_{i=1}^n\hat f(\{i\})\hat g(\{i\})\log\left(\frac{\|\frac12\partial_i f\|_2\|\frac12\partial_i g\|_2}
{\hat f(\{i\})\hat g(\{i\})}
\right)\ge-2W_1[f,g]\log\left(
\frac{\sum_i\|\frac12\partial_i f\|_2\|\frac12\partial_i g\|_2}{W_1[f,g]}\right).
$$
Parseval and Cauchy--Schwarz therefore give the following bound
\begin{equation}\label{eq:ordinary-derivative-bound}
W_3[f,g]\ge-\frac23W_1[f,g]\log\left(
\frac{\left(\sum_{S\ne\emptyset}|S|\hat f(S)^2\right)^{1/2}\left(\sum_{S\ne\emptyset}|S|\hat g(S)^2\right)^{1/2}}
     {W_1[f,g]}
\right).
\end{equation}
However, the quantities inside the logarithm are total derivative energies, not variances. They cannot simply be replaced by the smaller quantities $\Var(f)$ and $\Var(g)$, since that would increase the claimed lower bound.

To obtain variances, let $M$ be a real radial Fourier multiplier, i.e., $M\chi_S=m(|S|)\chi_S$. By \eqref{eq:derivative-fourier},
$$
\frac12\partial_iM(f-\E f)
=-\sum_{S\ni i}m(|S|)\hat f(S)\chi_{S\setminus\{i\}}.
$$
Consequently,
$$
\sum_{i=1}^n\left\|\frac12\partial_iM(f-\E f)\right\|_2^2=\sum_{i=1}^n\sum_{S\ni i}m(|S|)^2\hat f(S)^2=\sum_{S\ne\emptyset}|S|m(|S|)^2\hat f(S)^2=\sum_{k=1}^nk\,m(k)^2\sum_{|S|=k}\hat f(S)^2.
$$
For this expression to equal $\Var(f)$ for every real-valued $f$, each level must have coefficient $1$: taking $f=\chi_S$ with $|S|=k$ gives $k\,m(k)^2=1$. The positive choice is $m(k)=k^{-\frac{1}{2}}$, so $M=\N^{-1/2}$ on mean-zero functions. This choice also preserves the required first and third levels, because
$$
\begin{aligned}
&\sum_{i=1}^n\E\left[\frac12\partial_iM(f-\E f)\right]
                 \E\left[\frac12\partial_iM(g-\E g)\right]
=m(1)^2W_1[f,g],\\
&\sum_{i=1}^nW_2\left[
\frac12\partial_iM(f-\E f),\frac12\partial_iM(g-\E g)
\right]
=3m(3)^2W_3[f,g],
\end{aligned}
$$
and $m(1)^2=3m(3)^2=1$. Thus $\N^{-1/2}$ is the unique radial multiplier with positive eigenvalues that recovers the variance at every nonzero Fourier level.

\begin{defn}
For $f:\{0,1\}^n\to\mathbb R$ and $i\in[n]$, define the \emph{$i$-th discrete Riesz transform} by
\begin{equation}\label{eq:riesz-definition}
\mathfrak R_i f:=\frac12\partial_i\N^{-1/2}(f-\E f).
\end{equation}
\end{defn}

These are oriented versions of the discrete Riesz transforms associated with the Walsh number operator. Dimension-free norm estimates for these operators were studied by Lust-Piquard~\cite{LustPiquard1998,LustPiquard2004}, Junge, Mei and Parcet~\cite{JungeMeiParcet2018}, and Domelevo, Ivanisvili, Petermichl and Volberg~\cite{DomelevoEtAl2026}. Here we use their non-negativity on increasing functions and the identities in \cref{lem:fractional-identities}. 

\begin{lemma}\label{lem:fractional-positivity}
If $f$ is increasing, then $\mathfrak R_i f\ge0$ pointwise for every $i\in[n]$.
\end{lemma}

\begin{proof}
Using \eqref{eq:fractional-integral} on $f-\E f$, and then \eqref{eq:derivative-commutation}, we obtain
$$
\begin{aligned}
\mathfrak R_i f
&=\frac1{2\sqrt\pi}\int_0^\infty t^{-\frac{1}{2}}\partial_iP_t(f-\E f)\,dt=\frac1{2\sqrt\pi}\int_0^\infty t^{-\frac{1}{2}}\partial_iP_tf\,dt\\
&=\frac1{2\sqrt\pi}\int_0^\infty t^{-\frac{1}{2}}e^{-t}
P_t^{(-i)}(\partial_i f)\,dt\ge0.
\end{aligned}
$$
The derivative $\partial_i f$ is non-negative, $P_t^{(-i)}$ preserves non-negativity, and $t^{-\frac{1}{2}}e^{-t}>0$ for $t>0$.
\end{proof}

\begin{lemma}[Discrete Riesz transform identities]\label{lem:fractional-identities}
For $f,g:\{0,1\}^n\to\mathbb R$,
\begin{align}
\sum_{i=1}^n\E[\mathfrak R_i f]\,\E[\mathfrak R_i g]
&=W_1[f,g],\label{eq:identity-mean}\\
\sum_{i=1}^n\|\mathfrak R_i f\|_2^2
&=\Var(f),\qquad
\sum_{i=1}^n\|\mathfrak R_i g\|_2^2=\Var(g),\label{eq:identity-energy}\\
\sum_{i=1}^nW_2[\mathfrak R_i f,\mathfrak R_i g]
&=W_3[f,g].\label{eq:identity-level3}
\end{align}
\end{lemma}

\begin{proof}
By~\eqref{eq:derivative-fourier}, for $T\subseteq[n]\setminus\{i\}$,
\begin{equation}\label{eq:riesz-fourier}
\widehat{\mathfrak R_i f}(T)
=-\frac{\hat f(T\cup\{i\})}{\sqrt{|T|+1}},
\qquad
\widehat{\mathfrak R_i g}(T)
=-\frac{\hat g(T\cup\{i\})}{\sqrt{|T|+1}}.
\end{equation}
Taking $T=\emptyset$ gives $\E[\mathfrak R_i f]=-\hat f(\{i\})$ and $\E[\mathfrak R_i g]=-\hat g(\{i\})$. Their products sum to $W_1[f,g]$, proving \eqref{eq:identity-mean}. By Parseval,
$$
\sum_{i=1}^n\|\mathfrak R_i f\|_2^2=\sum_{i=1}^n\sum_{S\ni i}\frac{\hat f(S)^2}{|S|}=\sum_{S\ne\emptyset}\frac{\hat f(S)^2}{|S|}\sum_{i\in S}1=\sum_{S\ne\emptyset}\hat f(S)^2=\Var(f),
$$
and the identity for $g$ is identical. Finally, when $|T|=2$, the denominator in \eqref{eq:riesz-fourier} is $\sqrt3$. Hence
$$
\sum_{i=1}^nW_2[\mathfrak R_i f,\mathfrak R_i g]
=\frac13\sum_{i=1}^n
\sum_{\substack{S\ni i\\|S|=3}}
\hat f(S)\hat g(S)
=W_3[f,g],
$$
since each three-element set is counted three times.
\end{proof}

\begin{proof}[Proof of \cref{thm:level13-main}]
Assume first that $W_1[f,g]>0$. By \cref{lem:fractional-positivity}, every $\mathfrak R_i f$ and $\mathfrak R_i g$ is non-negative. Applying \eqref{eq:homogeneous-level-two-simple} coordinatewise and summing over all $i$ give
\begin{equation*}
    \begin{split}
\sum_{i=1}^nW_2[\mathfrak R_i f,\mathfrak R_i g]
&\ge -2\sum_{i=1}^n\E[\mathfrak R_i f]\,\E[\mathfrak R_i g]\log\left(
\frac{\|\mathfrak R_i f\|_2\|\mathfrak R_i g\|_2}{\E[\mathfrak R_i f]\,\E[\mathfrak R_i g]}\right)\\
&=-2W_1[f,g]\sum_{i=1}^n
\frac{\E[\mathfrak R_if]\,\E[\mathfrak R_i g]}{W_1[f,g]}\log\left(\frac{\|\mathfrak R_i f\|_2\|\mathfrak R_i g\|_2}
{\E[\mathfrak R_i f]\,\E[\mathfrak R_i g]}\right)\\
&\ge -2W_1[f,g]\log\left(\frac{\sum_i\|\mathfrak R_i f\|_2\|\mathfrak R_i g\|_2}{W_1[f,g]}\right)\\
&\ge -2W_1[f,g]\log\left(
\frac{\sqrt{\Var(f)\Var(g)}}{W_1[f,g]}
\right),
    \end{split}
\end{equation*}
where we use Jensen's inequality in the second inequality, and Cauchy--Schwarz together with~\cref{lem:fractional-identities} in the third inequality. Using~\eqref{eq:identity-level3}, we prove~\eqref{eq:level13-main} in this case.

If $W_1[f,g]=0$, then $\sum\limits_{i=1}^{n}\E[\mathfrak R_i f]\E[\mathfrak R_i g]=0.$
Since both functions are non-negative, for each $i$, either $\mathfrak R_i f=0$ or $\mathfrak R_i g=0$. Hence for every $i$, $W_2[\mathfrak R_i f,\mathfrak R_i g]=0$, and~\eqref{eq:identity-level3} gives $W_3[f,g]=0$.
\end{proof}

\begin{remark}
Let us compare with the Boolean consequence of Eldan's Gaussian Level-$1{:}3$ inequality. Let $f,g:\{0,1\}^n\to\mathbb R$ be increasing, and define
$$
F(x):=f(\mathbbm{1}_{\{x_1>0\}},\ldots,\mathbbm{1}_{\{x_n>0\}}),
\qquad
G(x):=g(\mathbbm{1}_{\{x_1>0\}},\ldots,\mathbbm{1}_{\{x_n>0\}}).
$$
Then $F,G$ are increasing and $\Var_{\gamma_n}(F)=\Var(f)$, with the analogous identity for $G$. For a standard Gaussian random variable $X$,
$$
\E[X\operatorname{sgn}(X)]=\sqrt{\frac{2}{\pi}},
\qquad
\E[(X^3-3X)\operatorname{sgn}(X)]=-\sqrt{\frac{2}{\pi}}.
$$
Since $\chi_{\{i\}}((\one_{\{x_j>0\}})_{j=1}^n)=-\operatorname{sgn}(x_i)$ almost everywhere, we obtain
$$
Q^{(1)}_i(F)=-\sqrt{\frac{2}{\pi}}\cdot\hat{f}(\{i\}),
\qquad
Q^{(3)}_{ijk}(F)=
\begin{cases}
-\left(\frac{2}{\pi}\right)^{\frac{3}{2}}\hat f(\{i,j,k\}),&i,j,k\text{ distinct},\\
0,&\text{exactly two indices coincide},\\
\sqrt{\frac{2}{\pi}}\hat{f}(\{i\}),&i=j=k.
\end{cases}
$$
Consequently,
\begin{equation}\label{eq:orthant-level13}
\begin{aligned}
\langle Q^{(1)}(F),Q^{(1)}(G)\rangle&=\frac{2}{\pi}W_1[f,g],\\
\langle Q^{(3)}(F),Q^{(3)}(G)\rangle&=6\left(\frac{2}{\pi}\right)^3W_3[f,g]+\frac{2}{\pi}W_1[f,g].
\end{aligned}
\end{equation}
Applying~\cite[Proposition~17]{Eldan2022}, in the form \eqref{eq:eldan-gaussian-level13}, therefore gives, for $A>0$,
\begin{equation}\label{eq:orthant-level13-corollary}
W_3[f,g]\ge-\frac{\pi^2}{24}W_1[f,g]
\left(1+e^8\log\left(\frac{\pi e\sqrt{\Var(f)\Var(g)}}{2W_1[f,g]}\right)\right).
\end{equation}
\end{remark}

\subsection{A two-function higher-level inequality}\label{subsec:keller-kindler}

For a bound at level $d$, we modify the Riesz transforms to obtain the factor $\frac{1}{d}$ in~\eqref{eq:keller-kindler-main}. The standard Riesz transforms preserve the variance, but the sum of the products of their $L^1$ norms is $W_1[f,g]$. Shifting the number operator by $d-1$ reduces this sum to $\frac{W_1[f,g]}{d}$ while keeping the total $L^2$ energy bounded by the variance and multiplying the target cross level by $\frac{d}{2d-1}$.

Fix an integer $d\ge2$ and increasing functions $f,g:\{0,1\}^n\to[-1,1]$. For $i\in[n]$, define
\begin{equation}\label{eq:shifted-riesz-definition}
\mathfrak R_{i,d}f:=\frac12\partial_i\bigl(\N+(d-1)\mathrm{Id}\bigr)^{-1/2}(f-\E f),
\end{equation}
and define $\mathfrak R_{i,d}g$ analogously. Here $\mathrm{Id}$ is the identity operator and
$$
\bigl(\N+(d-1)\mathrm{Id}\bigr)^{-1/2}\chi_S=(|S|+d-1)^{-1/2}\chi_S.
$$
As in \eqref{eq:fractional-integral} and \eqref{eq:derivative-commutation},
\begin{align}\label{eq:fractional-positivity}
\mathfrak R_{i,d}f=\frac1{2\sqrt\pi}\int_0^\infty t^{-1/2}e^{-dt}P_t^{(-i)}(\partial_i f)\,dt\ge0.
\end{align}
Since the semigroup preserves expectation, \eqref{eq:mean-derivative} gives
$$
\|\mathfrak R_{i,d}f\|_1=\E[\mathfrak R_{i,d}f]=\frac{\E[\partial_i f]}{2\sqrt d}=-\frac{\hat f(\{i\})}{\sqrt d},
$$
and similarly for $g$. Consequently,
\begin{equation}\label{eq:shifted-riesz-l1}
\sum_{i=1}^n\|\mathfrak R_{i,d}f\|_1\|\mathfrak R_{i,d}g\|_1=\frac{W_1[f,g]}d.
\end{equation}
For $T\subseteq[n]\setminus\{i\}$, \eqref{eq:derivative-fourier} gives
$$
\widehat{\mathfrak R_{i,d}f}(T)=-\frac{\hat f(T\cup\{i\})}{\sqrt{|T|+d}}.
$$
By Parseval,
$$
\sum_{i=1}^n\|\mathfrak R_{i,d}f\|_2^2=\sum_{S\ne\emptyset}\frac{|S|}{|S|+d-1}\hat f(S)^2\le\Var(f),
$$
and the same holds for $g$. Hence Cauchy--Schwarz and the boundedness of $f,g$ give
\begin{equation}\label{eq:shifted-riesz-l2}
\sum_{i=1}^n\|\mathfrak R_{i,d}f\|_2\|\mathfrak R_{i,d}g\|_2\le\sqrt{\Var(f)\Var(g)}\le1.
\end{equation}
Counting each set of size $d$ once for every coordinate it contains also gives
\begin{equation}\label{eq:shifted-riesz-level}
\sum_{i=1}^nW_{d-1}[\mathfrak R_{i,d}f,\mathfrak R_{i,d}g]=\frac{d}{2d-1}W_d[f,g].
\end{equation}

\begin{proof}[Proof of \cref{thm:keller-kindler-main}]
If $d>n$, then $W_d[f,g]=0$. If $W_1[f,g]=0$, \eqref{eq:shifted-riesz-l1} and \eqref{eq:fractional-positivity} imply that for each $i$ at least one of $\mathfrak R_{i,d}f,\mathfrak R_{i,d}g$ vanishes identically. Equation~\eqref{eq:shifted-riesz-level} then gives $W_d[f,g]=0$. We may therefore assume $W_1[f,g]>0$ and $2\le d\le n$.

Let $0<\rho<1$ and $\theta=\frac{1-\rho^2}{1+\rho^2}$. The standard Bonami--Beckner read-valued hypercontractivity~\cite{Beckner1975,Bonami1970} and $L^1$--$L^2$ interpolation give
$$
\|T_\rho\mathfrak R_{i,d}f\|_2\le\|\mathfrak R_{i,d}f\|_{1+\rho^2}
\le\|\mathfrak R_{i,d}f\|_1^\theta\|\mathfrak R_{i,d}f\|_2^{1-\theta},
$$
and similarly for $g$. Here $T_\rho$ acts on the cube without coordinate $i$. Multiplying the two estimates, summing over $i$, and applying H\"older's inequality, we obtain
$$
\begin{aligned}
\sum_{i=1}^n\|T_\rho\mathfrak R_{i,d}f\|_2\|T_\rho\mathfrak R_{i,d}g\|_2&\le\sum_{i=1}^n\bigl(\|\mathfrak R_{i,d}f\|_1\|\mathfrak R_{i,d}g\|_1\bigr)^\theta\bigl(\|\mathfrak R_{i,d}f\|_2\|\mathfrak R_{i,d}g\|_2\bigr)^{1-\theta}\\
&\le\left(\sum_{i=1}^n\|\mathfrak R_{i,d}f\|_1\|\mathfrak{R}_{i,d}g\|_1\right)^\theta\left(\sum_{i=1}^n\|\mathfrak R_{i,d}f\|_2\|\mathfrak R_{i,d}g\|_2\right)^{1-\theta}\\
&\le\left(\frac{W_1[f,g]}d\right)^\theta,
\end{aligned}
$$
where the last inequality uses~\eqref{eq:shifted-riesz-l1} and~\eqref{eq:shifted-riesz-l2}. Although the cross levels may be negative, \eqref{eq:shifted-riesz-level} and Cauchy--Schwarz at level $d-1$ give
$$
\begin{aligned}
\frac d{2d-1}|W_d[f,g]|&\le\sum_{i=1}^n\sqrt{W_{d-1}[\mathfrak R_{i,d}f]\,W_{d-1}[\mathfrak R_{i,d}g]}\le\rho^{-2(d-1)}\sum_{i=1}^n
\|T_\rho\mathfrak R_{i,d}f\|_2\|T_\rho\mathfrak R_{i,d}g\|_2\\
&\le\rho^{-2(d-1)}\left(\frac{W_1[f,g]}d\right)^\theta.
\end{aligned}
$$
Rearranging gives
\begin{equation}\label{eq:shifted-riesz-parameter-bound}
|W_d[f,g]|\le\frac{2d-1}{d^2}\,W_1[f,g]\,\rho^{-2(d-1)}\exp\left(\frac{2\rho^2}{1+\rho^2}\log\frac d{W_1[f,g]}\right).
\end{equation}
Taking the minimum of the right-hand side over $\rho$, we obtain approximately
$$
\rho^2=\frac{d-1}{2\log(d/W_1[f,g])}.
$$
The hypothesis $0< W_1[f,g] \le e^{-2(d-1)}$ ensures $0<\rho^2<\frac{1}{4}$, and
$$
\frac{2\rho^2}{1+\rho^2}\log\frac{d}{W_1[f,g]}\le d-1.
$$
Substituting into \eqref{eq:shifted-riesz-parameter-bound} proves \eqref{eq:keller-kindler-main}.
\end{proof}

\section{From level inequalities to correlation bounds}\label{sec:beta}
\begin{prop}[Averaged Level-$1{:}3$ inequality]\label{prop:averaged-level13}
For every $0\le p<1$ and any increasing real-valued function $f,g$,
\begin{equation}\label{eq:averaged-level13}
\E_{(J,z)}W_3[f_{J\to z},g_{J\to z}]
\ge-2\E_{(J,z)}W_1[f_{J\to z},g_{J\to z}]\log\left(\frac{\sqrt{\Var(f)\Var(g)}}{\E_{(J,z)}W_1[f_{J\to z},g_{J\to z}]}\right),
\end{equation}
where the right-hand side is interpreted as $0$ when the expected first level is zero.
\end{prop}

\begin{proof}
Apply \cref{thm:level13-main} to each restricted pair $(f_{J\to z},g_{J\to z})$. Suppose first that $\E_{(J,z)}W_1[f_{J\to z},g_{J\to z}]>0$. The quantities
$$
q_{J,z}:=
\frac{p^{|J|}(1-p)^{n-|J|}2^{-|J|}\,W_1[f_{J\to z},g_{J\to z}]}
     {\E_{(J,z)}W_1[f_{J\to z},g_{J\to z}]}
$$
are non-negative and sum to $1$ over $J\subseteq[n]$ and $z\in\{0,1\}^J$. Jensen's inequality gives
$$
\begin{aligned}
&\frac{\E_{(J,z)}\left[W_1[f_{J\to z},g_{J\to z}]\log\left(\frac{\sqrt{\Var(f_{J\to z})\Var(g_{J\to z})}}{W_1[f_{J\to z},g_{J\to z}]}\right)\right]}
{\E_{(J,z)}W_1[f_{J\to z},g_{J\to z}]}=\sum_{J,z}q_{J,z}\log\left(\frac{\sqrt{\Var(f_{J\to z})\Var(g_{J\to z})}}
{W_1[f_{J\to z},g_{J\to z}]}\right)\\
&\quad\le\log\left(\sum_{J,z}q_{J,z}
\frac{\sqrt{\Var(f_{J\to z})\Var(g_{J\to z})}}{W_1[f_{J\to z},g_{J\to z}]}\right)\le\log\left(
\frac{\E_{(J,z)}\sqrt{\Var(f_{J\to z})\Var(g_{J\to z})}}{\E_{(J,z)}W_1[f_{J\to z},g_{J\to z}]}
\right).
\end{aligned}
$$
In the last step, $W_1[f_{J\to z},g_{J\to z}]$ cancels against its factor in $q_{J,z}$; we then include the omitted non-negative variance terms. A term with zero first level contributes $0$ to the original expectation.

By Cauchy--Schwarz and \cref{lem:restriction-variance},
$$
\E_{(J,z)}\sqrt{\Var(f_{J\to z})\Var(g_{J\to z})}\le\sqrt{\E_{(J,z)}\Var(f_{J\to z})\,\E_{(J,z)}\Var(g_{J\to z})}\le\sqrt{\Var(f)\Var(g)}.
$$
Multiplying the Jensen estimate by the expected first level, and then averaging \eqref{eq:level13-main}, proves \eqref{eq:averaged-level13}.
\end{proof}

\begin{lemma}\label{lem:beta-extraction}
For every $0\le p<1$ and any increasing real-valued function $f,g$,
\begin{equation}\label{eq:beta-extraction}
h(\rho)=\rho W_1[f,g]+\rho^2W_2[f,g]+3\int_0^\rho\frac{(\rho-p)^2}{(1-p)^3}
\E_{(J,z)}W_3[f_{J\to z},g_{J\to z}]\,dp.
\end{equation}
Moreover,
\begin{equation}\label{eq:self-consistency}
h(\rho)\ge\rho W_1[f,g]+\rho^2W_2[f,g]-6\rho^2h(\rho)\log\left(\frac{\sqrt{\Var(f)\Var(g)}\log(1/(1-\rho))}{h(\rho)}
\right).
\end{equation}
The final term is interpreted as $0$ if $h(\rho)=0$. If $W_1[f,g]>0$, then $h(\rho)>0$ for every $0<\rho\le1$.
\end{lemma}

\begin{proof}
By \cref{lem:restriction-level},
$$
\E_{(J,z)}W_3[f_{J\to z},g_{J\to z}]
=(1-p)^3\sum_{k=3}^n\binom{k}{3}p^{k-3}W_k[f,g].
$$
For every $k\ge3$, expanding the square and integrating term by term gives
$$
\begin{aligned}
\int_0^\rho(\rho-p)^2p^{k-3}\,dp
&=\int_0^\rho\left(\rho^2p^{k-3}-2\rho p^{k-2}+p^{k-1}\right)dp\\
&=\rho^k\left(\frac1{k-2}-\frac2{k-1}+\frac1k\right)
=\frac{2\rho^k}{k(k-1)(k-2)}.
\end{aligned}
$$
Since $3\binom{k}{3}=\frac{k(k-1)(k-2)}{2}$, it follows that
$$
\begin{aligned}
3\int_0^\rho\frac{(\rho-p)^2}{(1-p)^3}
\E_{(J,z)}W_3[f_{J\to z},g_{J\to z}]\,dp&=3\sum_{k=3}^n\binom{k}{3}W_k[f,g]
\int_0^\rho(\rho-p)^2p^{k-3}\,dp\\
&=\sum_{k=3}^n\rho^kW_k[f,g]
=h(\rho)-\rho W_1[f,g]-\rho^2W_2[f,g].
\end{aligned}
$$
This proves \eqref{eq:beta-extraction}. 

To prove~\eqref{eq:self-consistency}, we may assume $h(\rho)>0$. Insert \eqref{eq:averaged-level13} into \eqref{eq:beta-extraction} and use
$\E_{(J,z)}W_1[f_{J\to z},g_{J\to z}]=(1-p)h'(p)$. We obtain
$$
h(\rho)\ge\rho W_1[f,g]+\rho^2W_2[f,g]-6\int_0^\rho\frac{(\rho-p)^2}{(1-p)^2}h'(p)
\log\left(\frac{\sqrt{\Var(f)\Var(g)}}{(1-p)h'(p)}\right)dp.
$$
The logarithm is non-negative wherever $h'(p)>0$: indeed, the first-level norm bound and \cref{lem:restriction-variance} give
$$
(1-p)h'(p)\le\E_{(J,z)}\sqrt{\Var(f_{J\to z})\Var(g_{J\to z})}
\le\sqrt{\Var(f)\Var(g)}.
$$
For $0\le p\le\rho<1$, we have $0\le\frac{\rho-p}{1-p}\le\rho$. Therefore the subtracted integral is at most
$$
6\rho^2\int_0^\rho h'(p)
\log\left(\frac{\sqrt{\Var(f)\Var(g)}}{(1-p)h'(p)}\right)dp.
$$

Since $h'\ge0$ and $\int_0^\rho h'(p)\,dp=h(\rho)$, the function $\frac{h'(p)}{h(\rho)}$ is a probability density on $[0,\rho]$. Jensen's inequality gives
$$
\begin{aligned}
\int_0^\rho h'(p)\log\left(\frac{\sqrt{\Var(f)\Var(g)}}{(1-p)h'(p)}\right)dp&\le h(\rho)\log\left(\int_{\{0<p<\rho:\,h'(p)>0\}}\frac{h'(p)}{h(\rho)}\frac{\sqrt{\Var(f)\Var(g)}}{(1-p)h'(p)}\,dp
\right)\\
&\le h(\rho)\log\left(\frac{\sqrt{\Var(f)\Var(g)}}{h(\rho)}\int_0^\rho\frac{dp}{1-p}\right)\\
&=h(\rho)\log\left(
\frac{\sqrt{\Var(f)\Var(g)}\log(1/(1-\rho))}{h(\rho)} \right).
\end{aligned}
$$
Substitution proves \eqref{eq:self-consistency}. 
\end{proof}

\begin{remark}[The Beta integral]
The termwise integration in the proof is the Beta integral. Recall that
$B(a,b)=\int_0^1u^{a-1}(1-u)^{b-1}\,du$ for $a,b>0$. For $k\ge3$, the substitution $p=\rho u$ gives
$$
\int_0^\rho(\rho-p)^2p^{k-3}\,dp=\rho^k\int_0^1u^{k-3}(1-u)^2\,du=\rho^kB(k-2,3)
=\rho^k\frac{\Gamma(k-2)\Gamma(3)}{\Gamma(k+1)}=\frac{\rho^k}{3\binom{k}{3}}.
$$
Thus the factor $3$ and the kernel $\frac{(\rho-p)^2}{(1-p)^3}$ in~\eqref{eq:beta-extraction} cancel the restriction multiplier $(1-p)^3\binom{k}{3}p^{k-3}$ after integration, leaving exactly $\rho^k$ for each $k\ge3$.
\end{remark}

We now prove~\cref{thm:antipodal-main} following the strategy in~\cite[Proof of Theorem~1.3]{KKM2016}.
\begin{proof}[Proof of \cref{thm:antipodal-main}]
If $\I[f,g]=0$, the right-hand side vanishes and the Harris--Kleitman inequality gives the result. Assume $\I[f,g]>0$, and choose
\begin{equation}\label{eq:delta-antipodal}
\delta
:=\frac1{4\sqrt{\log\left(
4e\sqrt{\Var(f)\Var(g)}/\I[f,g]
\right)}}.
\end{equation}
Since $W_1[f,g]=\frac{\I[f,g]}{4}\le\sqrt{\Var(f)\Var(g)}$, the logarithm in \eqref{eq:delta-antipodal} is at least $1$, and hence $\delta\le\frac{1}{4}$.

By antipodality, $W_2[f,g]=0$. We claim that
\begin{equation}\label{eq:H-two-thirds}
h(\delta)\ge\frac23\delta W_1[f,g].
\end{equation}
Suppose otherwise, and write
$$
h(\delta)=u\delta W_1[f,g]
\qquad\text{with }0<u<\frac23.
$$
Substituting $h(\delta)=u\delta W_1[f,g]$ in \eqref{eq:self-consistency} and dividing by $\delta W_1[f,g]>0$ gives
\begin{equation}\label{eq:u-ineq}
\begin{split}
1-u&\le6\delta^2u\log\left(\frac{\sqrt{\Var(f)\Var(g)}\log(1/(1-\delta))}{u\delta W_1[f,g]}
\right)\\
&\le6\delta^2u\Biggl(
\log\left(\frac{e\sqrt{\Var(f)\Var(g)}}{W_1[f,g]}\right)-1+\log\frac1{u(1-\delta)}\Biggr),
\end{split}
\end{equation}
where we use $\log(\frac{1}{1-\delta})=\int_0^\delta(1-p)^{-1}\,dp\le\frac{\delta}{1-\delta}$. Now let us consider, for $0<u\le2/3$, the function
$$
u\longmapsto u\left(
\log\frac{e\sqrt{\Var(f)\Var(g)}}{W_1[f,g]}-1
+\log\frac1{u(1-\delta)}
\right).
$$
Its first derivative is
$$
\log\frac{e\sqrt{\Var(f)\Var(g)}}{W_1[f,g]}-2+\log\frac1{u(1-\delta)},
$$
and its second derivative is $-\frac{1}{u}<0$. It extends continuously to $u=0$ with value $0$, since $u\log\frac{1}{u}\to0$. Its value is positive for $0<u\le\frac{2}{3}$, so its maximum on $[0,\frac{2}{3}]$ occurs at an interior critical point or at $u=\frac{2}{3}$.

At an interior critical point, the first derivative is zero, so the bracket in \eqref{eq:u-ineq} equals $1$. The function value there is $u\le\frac{2}{3}$, which is at most
$\frac23\log(e\sqrt{\Var(f)\Var(g)}/W_1[f,g])$ because this logarithm is at least $1$. At the right endpoint, $1-\delta\ge\frac{3}{4}$ gives
$$
\frac23\left(\log\frac{e\sqrt{\Var(f)\Var(g)}}{W_1[f,g]}-1+\log\frac3{2(1-\delta)}\right)\le\frac23\log\frac{e\sqrt{\Var(f)\Var(g)}}{W_1[f,g]}.
$$
Thus \eqref{eq:u-ineq} and \eqref{eq:delta-antipodal} imply
$$
\frac13<1-u\le6\delta^2\cdot\frac23
\log\frac{e\sqrt{\Var(f)\Var(g)}}{W_1[f,g]}=4\delta^2\log\frac{e\sqrt{\Var(f)\Var(g)}}{W_1[f,g]}=\frac14,
$$
a contradiction. This proves \eqref{eq:H-two-thirds}.

Finally, by \eqref{eq:cov-Hdelta} and $W_1[f,g]=\frac{\I[f,g]}{4}$,
$$
\Cov(f,g)\ge h(\delta)
\ge\frac23\delta W_1[f,g]=\frac1{24}\cdot\frac{\I[f,g]}{\sqrt{\log\left(4e\sqrt{\Var(f)\Var(g)}/\I[f,g]\right)}}.
$$
\end{proof}

\begin{cor}\label{cor:eldan-form}
Under the assumptions of \cref{thm:antipodal-main},
$$
\Cov(f,g)\ge\frac1{24}\cdot
\frac{\I[f,g]}{\sqrt{\log(2e/\I[f,g])}}.
$$
\end{cor}

\begin{proof}
Since $f,g$ are Boolean, $\sqrt{\Var(f)\Var(g)}\le1/4$. Hence
$$
\log\left(
\frac{4e\sqrt{\Var(f)\Var(g)}}{\I[f,g]}
\right)
\le\log\left(\frac{e}{\I[f,g]}\right)
\le\log\left(\frac{2e}{\I[f,g]}\right).
$$
Apply \cref{thm:antipodal-main}.
\end{proof}

For \cref{thm:second-order-main}, we retain the second-level term in \eqref{eq:self-consistency} and choose $\delta$ so that its negative part is at most $\delta W_1[f,g]/4$.

\begin{proof}[Proof of \cref{thm:second-order-main}]
If $W_1[f,g]=0$, the Harris--Kleitman inequality gives $\Cov(f,g)\ge0$. Assume $W_1[f,g]>0$, and choose
\begin{equation}\label{eq:delta-second-order}
\delta:=\frac14\min\left\{
\frac1{\sqrt{\log\left(e\sqrt{\Var(f)\Var(g)}/W_1[f,g]\right)}},\frac{W_1[f,g]}{\bigl(-W_2[f,g]\bigr)_+}\right\},
\end{equation}
with the usual convention when $W_2[f,g]\ge0$. Then $\delta\le1/4$ and
$\frac{\delta\bigl(-W_2[f,g]\bigr)_+}{W_1[f,g]}\le\frac14$. Write $h(\delta)=u\delta W_1[f,g]$, where $u>0$ by~\eqref{eq:beta-extraction} and the monotonicity of $h$. From \eqref{eq:self-consistency} and $\log\frac{1}{1-\delta}\le\frac{\delta}{1-\delta}$,
\begin{equation}\label{eq:second-u}
u\ge1-
\frac{\delta\bigl(-W_2[f,g]\bigr)_+}{W_1[f,g]}-6\delta^2u\Biggl(\log\left(\frac{e\sqrt{\Var(f)\Var(g)}}{W_1[f,g]}
\right)-1+\log\frac1{u(1-\delta)}
\Biggr).
\end{equation}
We claim that $u\ge\frac{1}{4}$. If $0<u<\frac{1}{4}$, use $u\log\frac{1}{u}\le\frac{1}{e}$, $\log\frac{1}{1-\delta}\le\log\frac{4}{3}$, and $$\log\left(\frac{e\sqrt{\Var(f)\Var(g)}}{W_1[f,g]}\right)\ge1$$ to obtain
$$
\begin{aligned}
&u\Biggl(\log\left(
\frac{e\sqrt{\Var(f)\Var(g)}}{W_1[f,g]}
\right)-1+\log\frac1{u(1-\delta)}\Biggr)\\
&\le\frac14\left[
\log\left(\frac{e\sqrt{\Var(f)\Var(g)}}{W_1[f,g]}\right)-1\right]+\frac1e+\frac14\log\frac43\le\log\left(\frac{e\sqrt{\Var(f)\Var(g)}}{W_1[f,g]}\right).
\end{aligned}
$$
Moreover, \eqref{eq:delta-second-order} implies
$$
\delta^2\le
\frac1{16\log\left(
 e\sqrt{\Var(f)\Var(g)}/W_1[f,g]
\right)}.
$$
Consequently, \eqref{eq:second-u} gives $u\ge1-\frac{1}{4}-\frac{3}{8}=\frac{3}{8}$, a contradiction. Hence $h(\delta)\ge\delta W_1[f,g]/4$. Combining this with \eqref{eq:cov-Hdelta} and \eqref{eq:delta-second-order} proves \eqref{eq:second-order-main}.
\end{proof}

\section{Sharpness and concluding remarks}\label{sec:sharpness}

\begin{proof}[Proof of \cref{prop:sharpness}]
Let $n$ be odd. Define
$$
f(x):=\mathbbm{1}_{\{x_1+\cdots+x_n\ge1\}},
\qquad
g(x):=\mathbbm{1}_{\{x_1+\cdots+x_n\ge(n+1)/2\}}.
$$
Thus $f$ is the OR function and $g$ is majority. Both are increasing, and $g$ is antipodal. Since $g=1$ implies $f=1$,
$$
\Cov(f,g)
=\mathbb E[g]-\mathbb E[f]\mathbb E[g]
=\frac12-\left(1-2^{-n}\right)\frac12
=2^{-n-1}.
$$
For every coordinate,
$$
\Inf_i[f]=2^{-(n-1)},
\qquad
\Inf_i[g]
=2^{-(n-1)}\binom{n-1}{(n-1)/2}.
$$
Consequently,
$$
\I[f,g]
=n\,2^{-2(n-1)}\binom{n-1}{(n-1)/2}
\asymp\sqrt n\,2^{-n},
$$
where we used the standard central-binomial estimate. Hence $\log\left(\frac{2e}{\I[f,g]}\right)\asymp n$. For every $\alpha<1/2$,
$$
\frac{\Cov(f,g)\log^\alpha(2e/\I[f,g])}{\I[f,g]}
\asymp n^{\alpha-1/2}\longrightarrow0.
$$
Thus no positive universal constant can make the asserted bound hold.
\end{proof}

The constant $2$ in \cref{thm:level13-main} is not expected to be optimal. To make the question precise, define
$$
C_*:=\sup
\frac{\bigl(-W_3[f,g]\bigr)_+}
{W_1[f,g]\log\left(\sqrt{\Var(f)\Var(g)}/W_1[f,g]\right)},
$$
where the supremum ranges over increasing real-valued pairs with a positive denominator. Our theorem gives $C_*\le2$. A three-variable example gives a nontrivial lower bound. Indeed, for $0\le\lambda\le1$, let
$$
f_\lambda=-\chi_{\{1\}}-\chi_{\{2\}}-\chi_{\{3\}}-\lambda\chi_{\{1,2,3\}},
$$
$$
g_\lambda=-\chi_{\{1\}}-\chi_{\{2\}}-\chi_{\{3\}}+\lambda\chi_{\{1,2,3\}}.
$$
Both functions are increasing, and
$$
W_1[f_\lambda,g_\lambda]=3,
\qquad
W_3[f_\lambda,g_\lambda]=-\lambda^2,
\qquad
\Var(f_\lambda)=\Var(g_\lambda)=3+\lambda^2.
$$
Taking $\lambda=1$ yields
\begin{equation}\label{eq:Cstar-lower}
C_*\ge\frac1{3\log(4/3)}=1.158686\ldots.
\end{equation}
We ask whether the three-variable example determines $C_*$.

\begin{ques}\label{ques:optimal-level13}
Is $C_*=1/(3\log(4/3))$? Equivalently, does every increasing real-valued functions satisfy
$$
W_3[f,g]
\ge-\frac{W_1[f,g]}{3\log(4/3)}
\log\left(
\frac{\sqrt{\Var(f)\Var(g)}}{W_1[f,g]}
\right)?
$$
\end{ques}

The two steps in the proof of \cref{thm:level13-main} are the imaginary-noise estimate of \cref{lem:imaginary-level-two} and its application to the non-negative discrete Riesz transforms. The same Fourier calculation as in \cref{lem:fractional-identities} gives
$$
\sum_{i=1}^nW_{k-1}[\mathfrak R_i f,\mathfrak R_i g]=W_k[f,g]
\qquad(1\le k\le n).
$$
A possible extension is to combine this identity with higher-level estimates for non-negative functions to obtain Level-$1{:}(2r+1)$ inequalities.

\section*{Acknowledgements.}

This project began during the spring semester of 2025, when the first
author visited Shandong University and participated in the seminar \emph{Topics in Combinatorics / Analysis of Boolean Functions}. The project had two main motivations. The first was to investigate whether Eldan's methods~\cite{Eldan2022} could be used to prove the correlation inequality equivalent to Chv\'atal's conjecture.
The second was to seek a more direct proof of Eldan's correlation
inequalities, inspired by subsequent proofs of Talagrand's isoperimetric
conjecture~\cite{Talagrand1997} that avoid the pathwise stochastic
analysis used in Eldan and Gross's original proof~\cite{EG2022};
see~\cite{EKLM2025,IvanisviliZhang2026,Rosenthal2020}.

The authors acknowledge the use of AI tools. At a later stage of the
project, they used ChatGPT to help understand Eldan's proof~\cite{Eldan2022} and identify steps that might be replaced.
In particular, ChatGPT 5.6 suggested using Beckner's complex hypercontractivity in place of Eldan's entropy estimate. The ideas of using discrete Riesz transforms and random restrictions originated entirely with the authors, who developed all the proofs.
ChatGPT Astra 6.0 was used for proofreading and to discuss the
optimal constant in~\cref{thm:level13-main}. The authors take full responsibility for the content of the paper.

\end{document}